\documentclass{amsart}
\usepackage{amsmath,amssymb}
\usepackage[dvips]{graphics}
\usepackage[all]{xy}
\newtheorem{Theorem}{Theorem}[section]
\newtheorem{Lemma}[Theorem]{Lemma}
\newtheorem{Proposition}[Theorem]{Proposition}

\newtheorem{Example}[Theorem]{Example}

\newtheorem{Remark}[Theorem]{Remark}

\makeatletter
\@addtoreset{figure}{section}
\def\@thmcountersep{-}
\makeatother

\numberwithin{equation}{section}

\begin{document}

\title{Mapping a knot by a continuous map}

\author{Kouki Taniyama}
\address{Department of Mathematics, School of Education, Waseda University, Nishi-Waseda 1-6-1, Shinjuku-ku, Tokyo, 169-8050, Japan}
\email{taniyama@waseda.jp}
\thanks{The author was partially supported by Grant-in-Aid for Scientific Research (C) (No. 24540100), Japan Society for the Promotion of Science.}


\date{}

\dedicatory{}


\begin{abstract}
By a fixed continuous map from a $3$-space to itself, a knot in the $3$-space may be mapped to another knot in the $3$-space. 
We analyze possible knot types of them. Then we map a knot repeatedly by a fixed continuous map and analyze possible infinite sequences of knot types. 
\end{abstract}

\maketitle

\section{Introduction} 
%

Let ${\mathbb R}^3$ be the $3$-dimensional Euclidean space. Let $f:{\mathbb R}^3\to{\mathbb R}^3$ be a continuous map. Let $k$ be an oriented knot in ${\mathbb R}^3$. Then $k$ is a subset of ${\mathbb R}^3$ that is homeomorphic to a circle. Suppose that the restriction map $f|_k:k\to{\mathbb R}^3$ is injective. Then the image $f(k)$ is again an oriented knot in ${\mathbb R}^3$. 
In this paper we consider the problem whether or not the knot types of $k$ and $f(k)$ can be chosen arbitrary. We also consider the knot types of the infinite sequence of oriented knots $k,f(k),f^2(k),\cdots$. 

A knot in ${\mathbb R}^3$ is usually mapped by an orientation preserving self-homeomorphism of ${\mathbb R}^3$. Mapping a knot by a non-injective continuous map is unusual and seems to be a mismatch. We are simply interested in seeing what happens when we relax the condition of homeomorphism to continuous map. Most of the results in this paper states the existence of various examples. After all we see that almost everything can happen when we map a knot by a continuous map. 
We hope there will be further studies and/or applications of our study. 

\vskip 3mm

\section{Folding a knot}\label{section2} 

We begin with the following imaginary stories. Suppose that there is a right circle drawn on a transparent paper. Fold the paper along a straight line intersecting the circle at two points but missing the center of the circle. Then we find a shape which is no more a right circle but yet a topological circle. Next suppose that a person living in a $4$-dimensional space draw a knot on a transparent paper which is of course $3$-dimensional, and fold the paper along a flat plane. Then he/she find a new knot. Here we have a question what are the knot types before and after folding the paper. We formulate and answer the question as follows. 

We denote the set of all oriented tame knot types in the $3$-dimensional Euclidean space ${\mathbb R}^3$ by ${\mathcal K}$. Then an element $K$ of ${\mathcal K}$ is an oriented knot type in ${\mathbb R}^3$. Note that $K$ is an ambient isotopy class of oriented knots in ${\mathbb R}^3$. Therefore an element $k$ of $K$ is an oriented knot in ${\mathbb R}^3$. 
Since $K$ is a tame knot type there is an element $k$ of $K$ that is smooth or polygonal. Here a knot in ${\mathbb R}^3$ is {\it polygonal} if it is a union of finitely many straight line segments. All oriented knots in this paper are smooth or polygonal unless otherwise stated. 
Let $F:{\mathbb R}^3\to{\mathbb R}^3$ be a map defined by $F(x,y,z)=(x,y,|z|)$.

\vskip 3mm

\begin{Theorem}\label{folding a knot}
Let $K_1$ and $K_2$ be elements of ${\mathcal K}$. 
Then there is an element $k_1$ of $K_1$ such that $F$ maps $k_1$ homeomorphically onto an element $F(k_1)$ of $K_2$. 
\end{Theorem}

\vskip 3mm

\begin{figure}[htbp]
      \begin{center}
\scalebox{0.6}{\includegraphics*{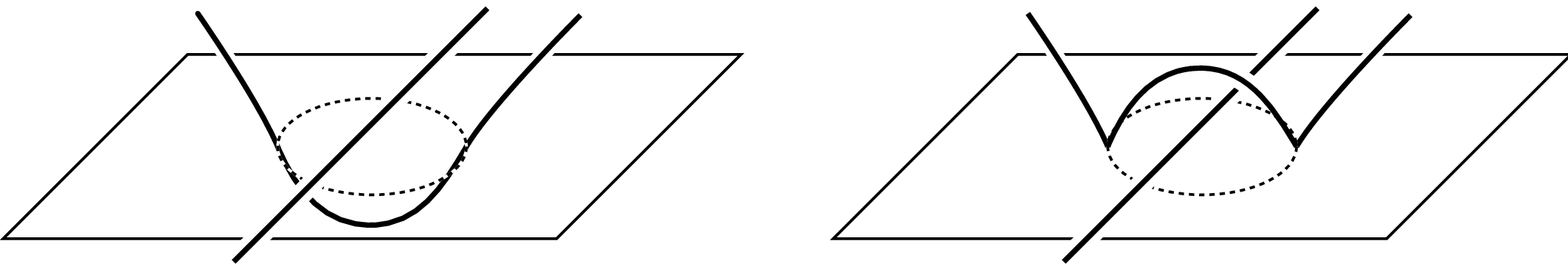}}
      \end{center}
   \caption{}
  \label{under-arc}
\end{figure}

\begin{figure}[htbp]
      \begin{center}
\scalebox{0.6}{\includegraphics*{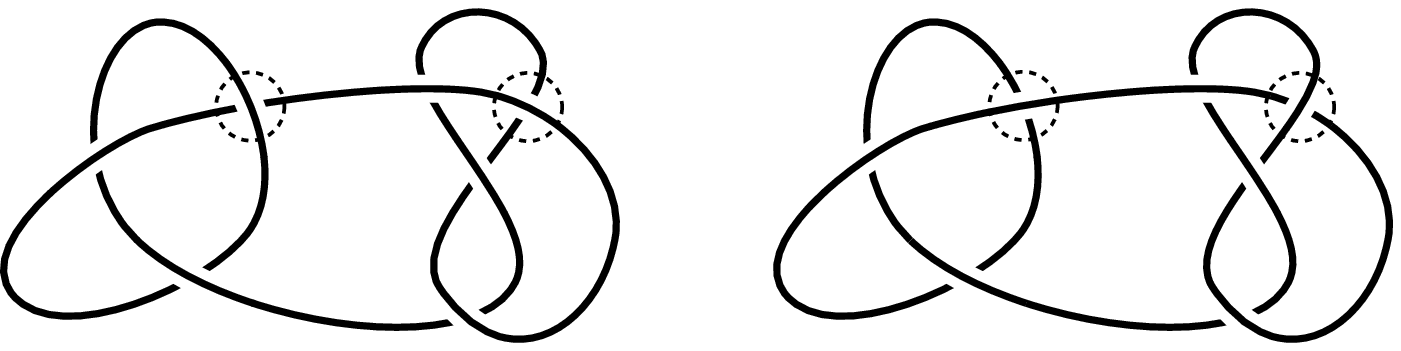}}
      \end{center}
   \caption{}
  \label{folding}
\end{figure}

\noindent{\bf Proof.} Let $D_i$ be an oriented knot diagram representing $K_i$ and ${\mathcal C}_i$ a set of crossings of $D_i$ such that changing over/under at all crossings in ${\mathcal C}_i$ will turn $D_i$ into a trivial knot diagram $D_i'$ for $i=1,2$. Let ${\mathcal C}_i'$ be a set of crossings of $D_i'$ corresponding to ${\mathcal C}_i$ for $i=1,2$. 
Let $D$ be a knot diagram obtained by a diagram-connected sum of $D_1$ and $D_2'$. 
Let $k$ be an oriented knot whose diagram is $D$. 
Then $k$ is an element of $K_1$. Suppose that $k$ is slightly above the $xy$-plane except the arcs corresponding to the under arcs of the crossings in ${\mathcal C}_1$ and ${\mathcal C}_2'$ so that the diagram of $F(k)$ is obtained from $D$ by changing over/under at each crossing in ${\mathcal C}_1$ and ${\mathcal C}_2'$. See Figure \ref{under-arc}. 
Then the diagram of $F(k)$ is a diagram-connected sum of $D_1'$ and $D_2$. Therefore $F(k)$ is an element of $K_2$ as desired. See for example Figure \ref{folding}. 
$\Box$

\section{generalization}\label{generalization} 

We now generalize Theorem \ref{folding a knot} as follows. By a {\it simple arc} (resp. {\it disk}, {\it 3-ball}) we mean a topological space homeomorphic to a closed interval (resp. closed 2-disk, closed 3-ball). 

\vskip 3mm

\begin{Theorem}\label{continuous map}
Let $f:{\mathbb R}^3\to{\mathbb R}^3$ be a continuous map. Suppose that there exist a simple arc $A$ and a 3-ball $B$ in ${\mathbb R}^3$ with the following properties. 

(1) $A\cap B=\partial A\cap\partial B=\partial A$,

(2) $A\cup A'$ is a trivial knot where $A'$ is a simple arc in $\partial B$ with $\partial A'=\partial A$,

(3) $f$ maps each of $A$ and $B$ homeomorphically onto its image,

(4) there is a simple arc $A''$ in ${\rm int}A$ such that $f(A)\cap f(B)=f(\partial A)\cup f(A'')$ and $(f(B),f(A''))$ is a trivial ball-arc pair,

(5) there are mutually disjoint disks $D_1$ and $D_2$ in ${\mathbb R}^3\setminus f({\rm int}B)$ such that $(D_1\cup D_2)\cap f(B)=(\partial D_1\cup\partial D_2)\cap f(\partial B)$ is a disjoint union of two simple arcs and $(\partial D_1\cup\partial D_2)\setminus f(B)=f(A)\setminus f(B)$. 

Let $K_1$ and $K_2$ be elements of ${\mathcal K}$. 
Then there is an element $k_1$ of $K_1$ such that $f$ maps $k_1$ homeomorphically onto an element $f(k_1)$ of $K_2$. 
\end{Theorem}

\vskip 3mm

We use the following lemma for the proof of Theorem \ref{continuous map}. This lemma is a version of Terasaka-Suzuki lemma. 
A first version of this lemma appeared in an expository paper \cite{Terasaka} written in Japanese. Then it is shown in \cite[Lemma 1]{Suzuki}. 
See also \cite[Lemma 2.1]{Yamamoto} and \cite[Lemma 2.1(1)]{Taniyama-Yasuhara}. 
For knot types $K_1$ and $K_2$, $d_G(K_1,K_2)$ denotes the Gordian distance between them. Namely $d_G(K_1,K_2)$ is the minimal number of crossing changes that is needed to deform $K_1$ to $K_2$.

\vskip 3mm

\begin{figure}[htbp]
      \begin{center}
\scalebox{0.6}{\includegraphics*{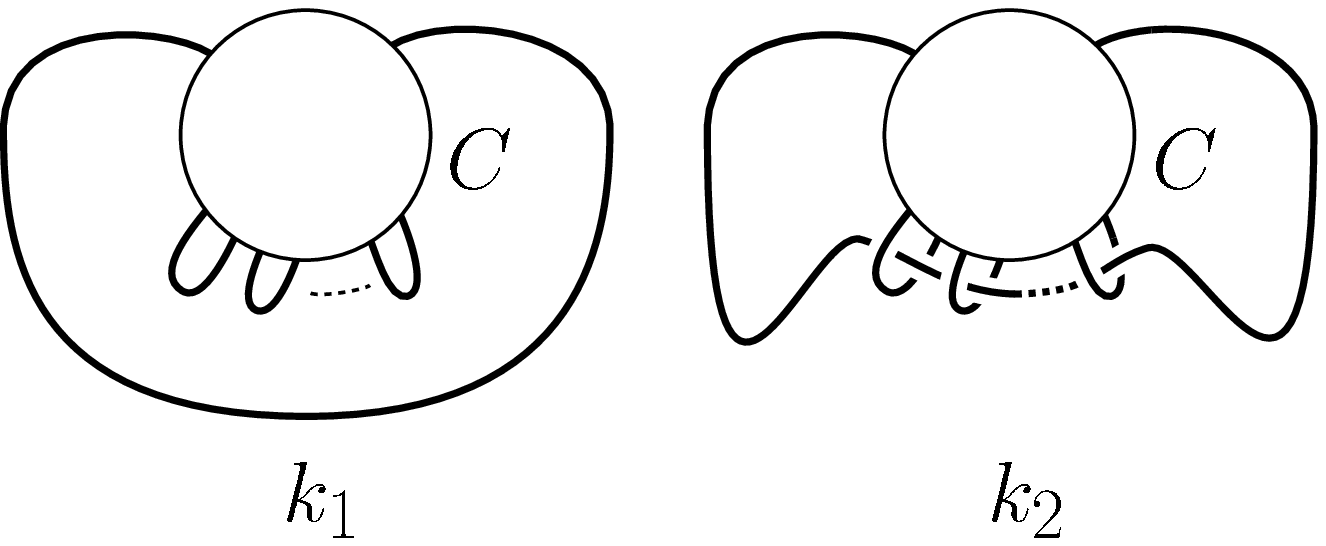}}
      \end{center}
   \caption{}
  \label{Terasaka-Suzuki}
\end{figure}

\begin{Lemma}\label{TS}
Let $K_1$ and $K_2$ be elements of ${\mathcal K}$. 
Then there is an element $k_1$ of $K_1$ and an element $k_2$ of $K_2$ with the following properties. 

(1) There is a $3$-ball $C$ in ${\mathbb R}^3$ such that $C\cap k_1=C\cap k_2$ and the pair $(C,C\cap k_1)$ is a $(d_G(K_1,K_2)+1)$-string tangle,

(2) Outside of $C$ is exactly as illustrated in Figure \ref{Terasaka-Suzuki}.

\end{Lemma}

\vskip 3mm

\begin{figure}[htbp]
      \begin{center}
\scalebox{0.6}{\includegraphics*{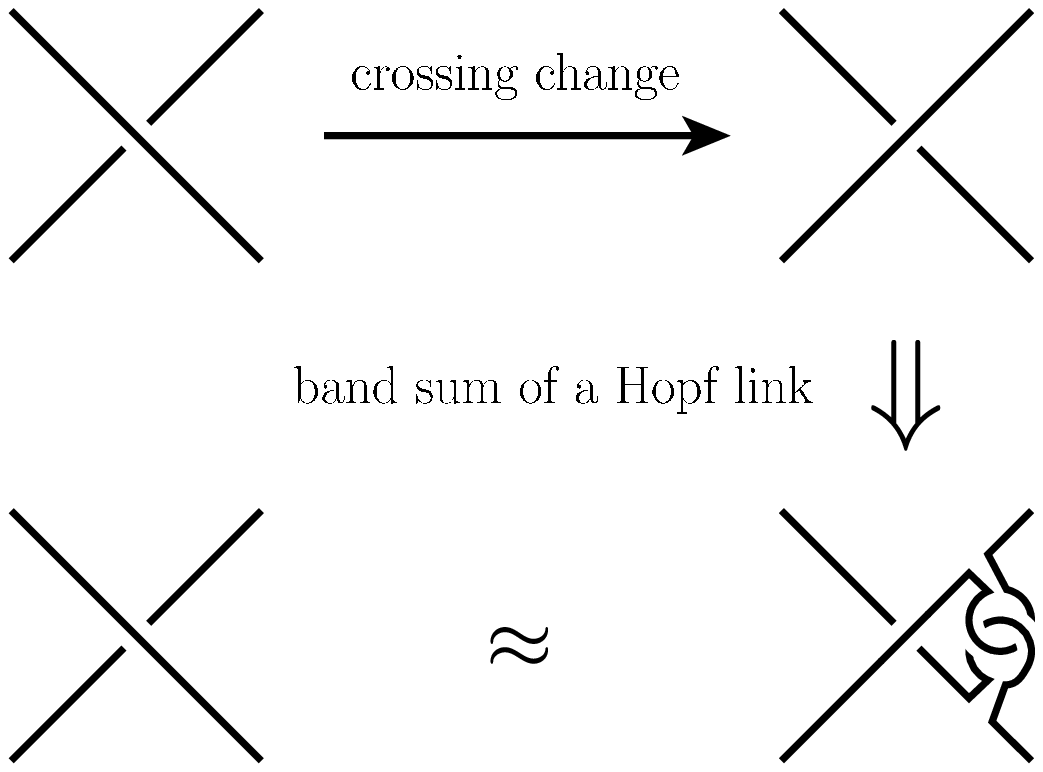}}
      \end{center}
   \caption{}
  \label{band-sum}
\end{figure}

\begin{figure}[htbp]
      \begin{center}
\scalebox{0.6}{\includegraphics*{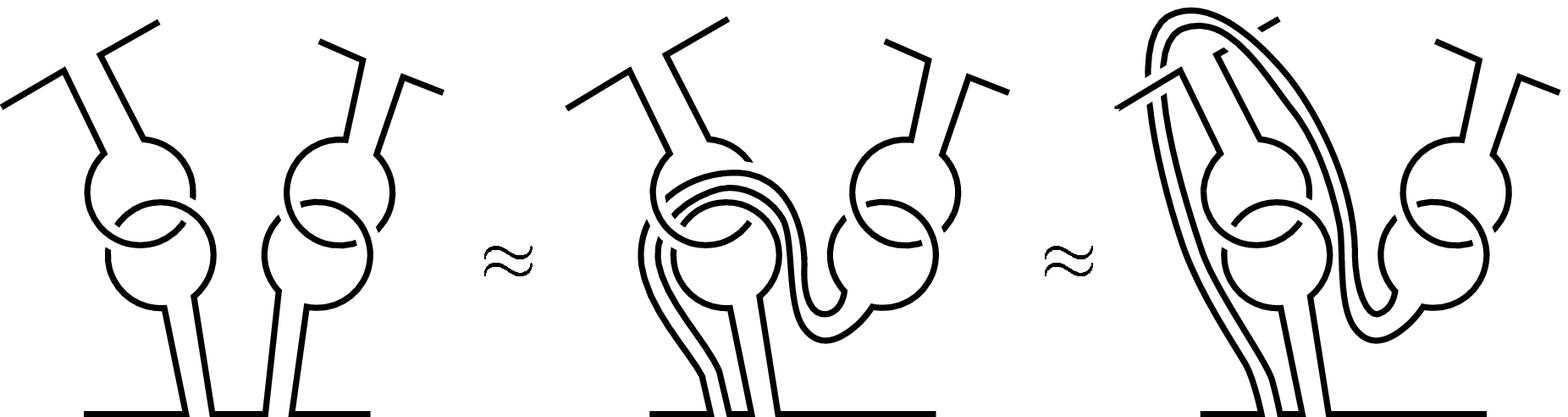}}
      \end{center}
   \caption{}
  \label{band-deformation}
\end{figure}

\begin{figure}[htbp]
      \begin{center}
\scalebox{0.6}{\includegraphics*{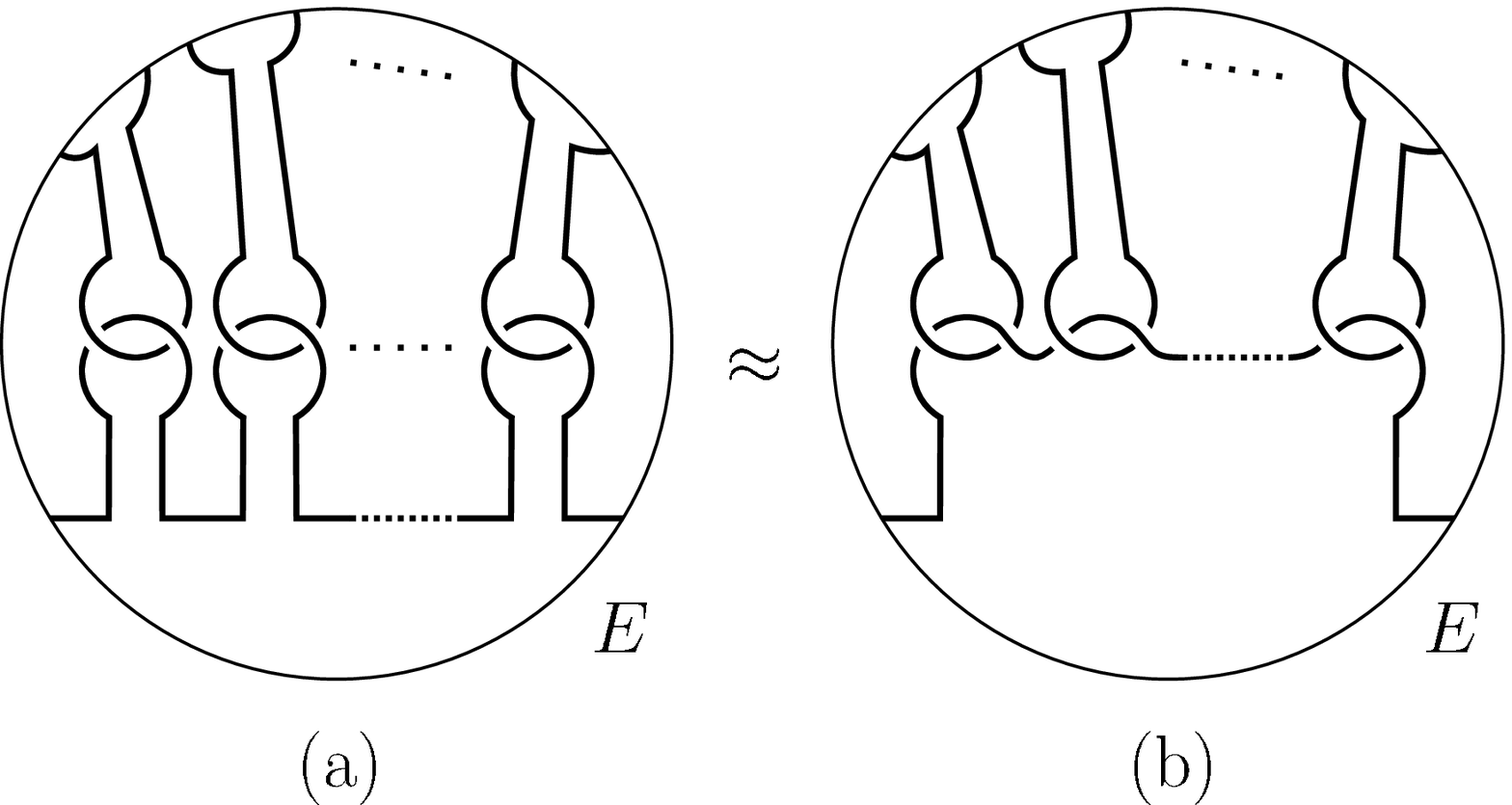}}
      \end{center}
   \caption{}
  \label{Hopf-links}
\end{figure}

\noindent{\bf Proof.} We consider a sequence of $d_G(K_2,K_1)=d_G(K_1,K_2)$ times crossing changes from $K_2$ to $K_1$. 
We replace a crossing change by a band sum of a Hopf link as illustrated in Figure \ref{band-sum}. 
Repeating the replacements we have an element $k_2$ of $K_2$ that is a band-sum of $d_G(K_1,K_2)$ Hopf links 
and an element $k_1$ of $K_1$. See \cite[Lemma 2.1(1)]{Taniyama-Yasuhara} for more detail. By using the deformation illustrated in Figure \ref{band-deformation} if necessary we can deform the band-sum without changing the knot type so that all Hopf links are contained in a $3$-ball, say $E$, where they are exactly as illustrated in Figure \ref{Hopf-links} (a). After deforming it as illustrated in Figure \ref{Hopf-links} (b) by an ambient isotopy, we see that it is equal to the outside of $C$ of $k_2$ in Figure \ref{Terasaka-Suzuki} in the $3$-sphere that is the one point compactification of ${\mathbb R}^3$. Therefore the outside of $E$ in the $3$-sphere give rise to the desired $(d_G(K_1,K_2)+1)$-string tangle. 
$\Box$


\vskip 3mm

\begin{figure}[htbp]
      \begin{center}
\scalebox{0.6}{\includegraphics*{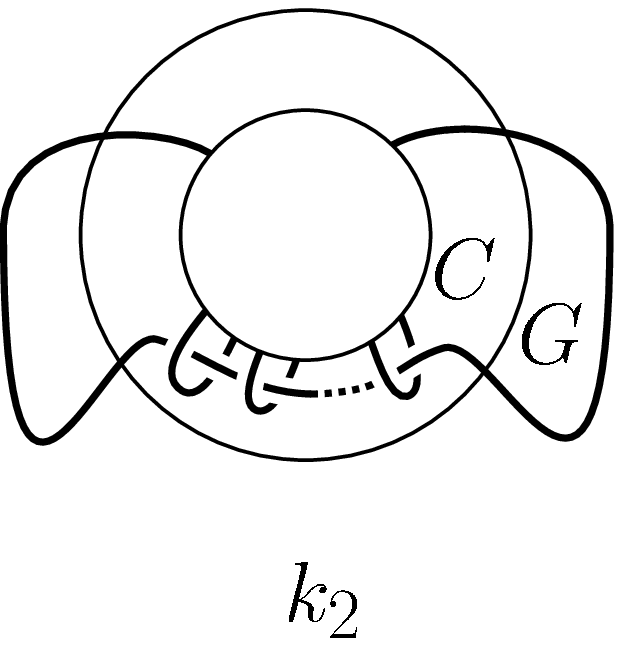}}
      \end{center}
   \caption{}
  \label{Terasaka-Suzuki2}
\end{figure}

\noindent{\bf Proof of Theorem \ref{continuous map}.} First suppose that the restriction map $f|_B$ of $f$ on $B$ preserves the orientation. 
We use $k_2$ in Lemma \ref{TS}. 
Let $G$ be a $3$-ball in ${\mathbb R}^3$ containing $C$ as illustrated in Figure \ref{Terasaka-Suzuki2}. 
Let $\alpha$ be the component of $k_2\cap G$ that is disjoint from $C$. 
By the assumption there exists an orientation preserving self-homemorphism $h:{\mathbb R}^3\to{\mathbb R}^3$ such that $h(G\cup k_2)=f(A\cup B)$ and $h(\alpha)=f(A'')$. 
Let $\beta=B\cap f^{-1}(h(G\cap k_2)\setminus\alpha)$ and $k_1=A\cup\beta$. Then we see that $k_1$ is an element of $K_1$ and $f(k_1)=h(k_2)$ is an element of $K_2$ as desired. 
Next suppose that $f|_B$ reverses the orientation. In this case we apply Lemma \ref{TS} for the mirror image of $K_1$ instead of $K_1$. Then the same proof above works. 
$\Box$

\vskip 3mm

\begin{figure}[htbp]
      \begin{center}
\scalebox{0.6}{\includegraphics*{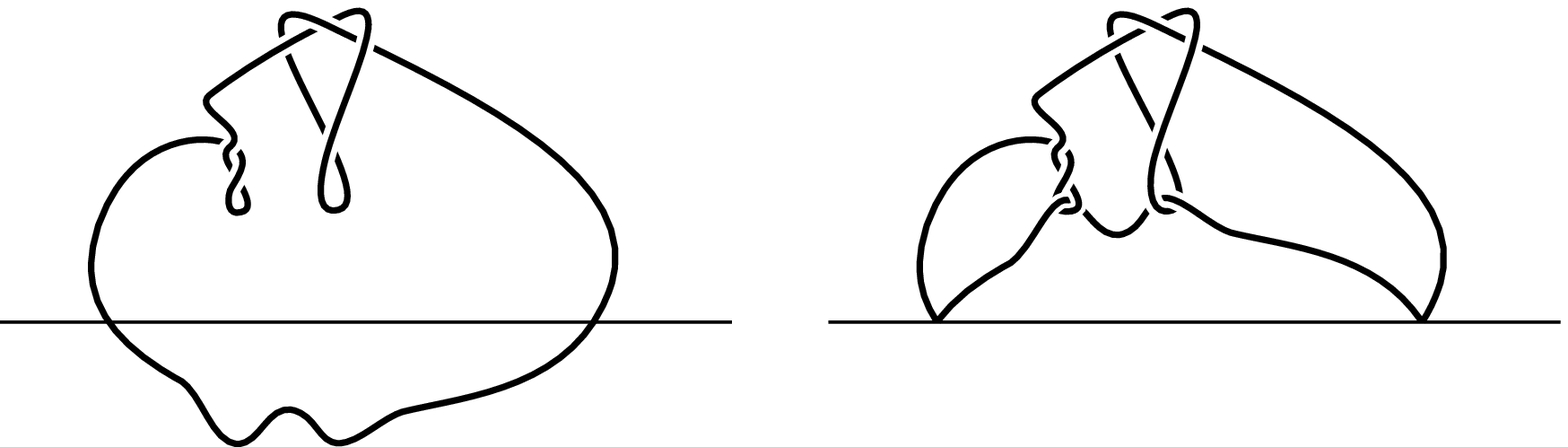}}
      \end{center}
   \caption{}
  \label{folding-example}
\end{figure}

It is easy to see that the map $F$ in section \ref{section2} satisfies the conditions of Theorem \ref{continuous map}. Therefore Theorem \ref{folding a knot} follows also from Theorem \ref{continuous map}. 
An example of $k$ and $F(k)$ along the proof of Theorem \ref{continuous map} is illustrated in Figure \ref{folding-example}. Here $K_1$ is the right-handed trefoil knot type and $K_2$ is the figure eight knot type.

\begin{figure}[htbp]
      \begin{center}
\scalebox{0.4}{\includegraphics*{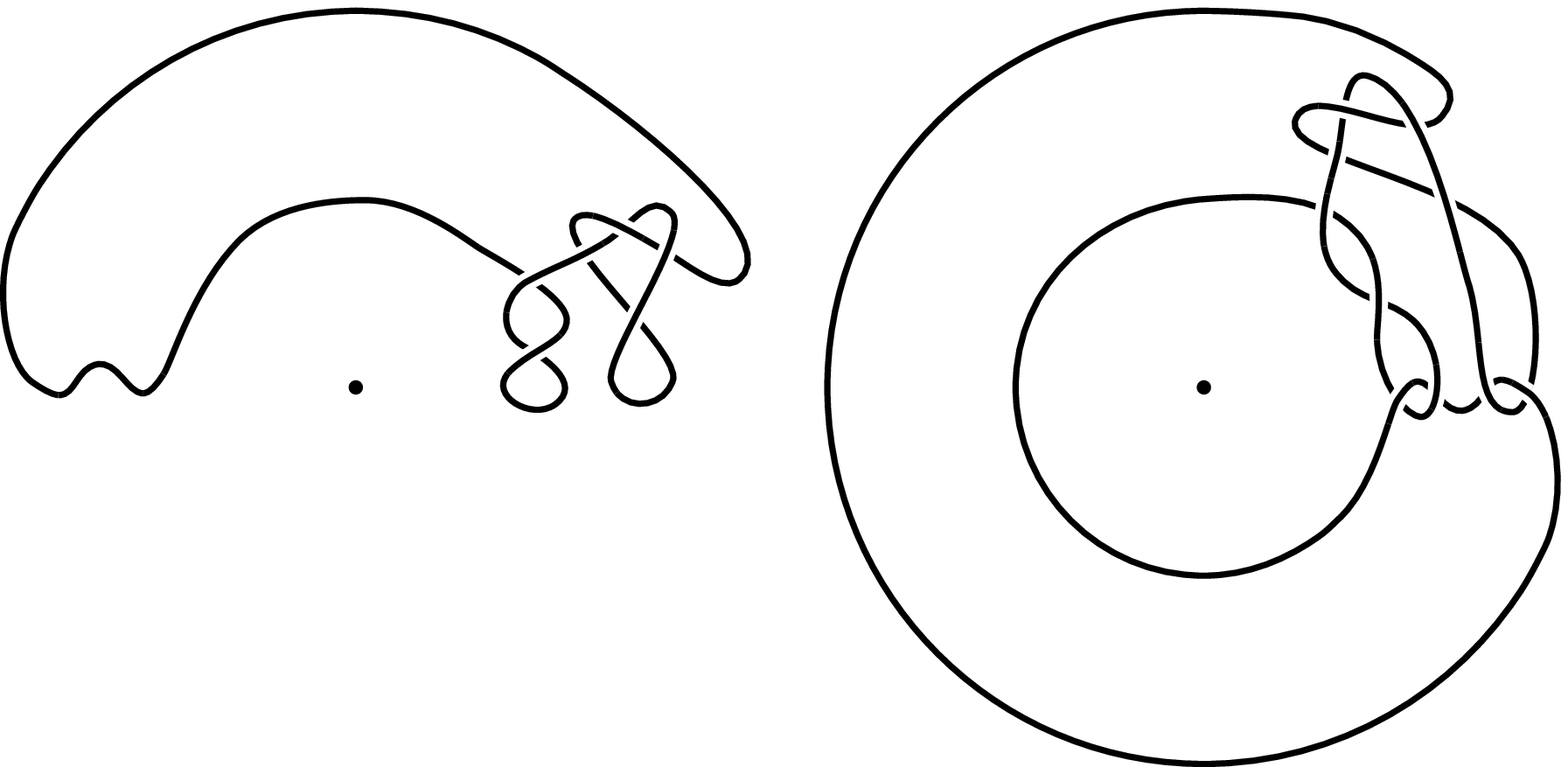}}
      \end{center}
   \caption{}
  \label{winding-example}
\end{figure}

As another example of Theorem \ref{continuous map} we consider the following continuous map. 
Let $W:{\mathbb R}^3\to{\mathbb R}^3$ be a map defined by $W(r{\rm cos}\theta,r{\rm sin}\theta,z)=(r{\rm cos}2\theta,r{\rm sin}2\theta,z)$. 
It is easy to see that the map $W$ also satisfies the conditions of Theorem \ref{continuous map}. 
An example is illustrated in Figure \ref{winding-example}. Here $K_1$ is the right-handed trefoil knot type and $K_2$ is the figure eight knot type.

\vskip 3mm

\begin{Remark}\label{remark}
{\rm 
The point of Theorem \ref{continuous map} is that a continuous map $f:{\mathbb R}^3\to{\mathbb R}^3$ is previously fixed. We will show in Proposition \ref{neighbourhood} that for given two oriented knots $k_1$ and $k_2$ it is easy to construct a continuous map $f:{\mathbb R}^3\to{\mathbb R}^3$ with $f(k_1)=k_2$. 
}
\end{Remark}

\vskip 3mm

\section{Mapping a knot repeatedly by a continuous map} 

Let $X$ be a set and $f:X\to X$ a map. We define a map $f^n:X\to X$ inductively by $f^0={\rm id}_{X}$, $f^1=f$ and $f^n=f\circ f^{n-1}$ for each $n\geq2$. Here $f\circ g$ denotes the composition map of $g$ and $f$. 
Let $A$ be a subset of $X$. We say that $A$ is {\it iteratively injective with respect to $f$} if the restriction map $f^n|_A:A\to X$ is injective for every non-negative integer $n$. 
Let $f:{\mathbb R}^3\to{\mathbb R}^3$ be a continuous map. Let $k$ be an oriented knot in ${\mathbb R}^3$ that is iteratively injective with respect to $f$. 
Then $f^n(k)$ is an oriented knot for each $n$. We analyze the knot types of them. We begin with the following example.

\vskip 3mm

\begin{Example}\label{torus knot}
{\rm
Let $k=\{((2+{\rm cos}3\theta){\rm cos}\theta,(2+{\rm cos}3\theta){\rm sin}\theta,{\rm sin}3\theta)|0\leq\theta\leq2\pi\}$ be a standard $(1,3)$-torus knot. Then $W^n(k)$ is a standard $(2^n,3)$-torus knot. Thus we have an infinite sequence $k,W(k),W^2(k),\cdots$ of knots such that no two of them belong to the same knot type. 
}
\end{Example}

\vskip 3mm

Now we state the following proposition for further examples. 

\vskip 3mm

\begin{Proposition}\label{neighbourhood}
Let $k_1$ and $k_2$ be oriented knots in ${\mathbb R}^3$. Let $V_1$ and $V_2$ be their tubular neighbourhoods in ${\mathbb R}^3$ respectively. Let $B$ be a $3$-ball such that $V_1\cup V_2$ is contained in ${\rm int}B$. Then there exists a continuous map $f:{\mathbb R}^3\to{\mathbb R}^3$ with the following properties. 

\begin{enumerate}
\item[(1)] $f$ maps the pair $(V_1,k_1)$ homeomorphically onto the pair $(V_2,k_2)$, 
\item[(2)] $f({\bf x})={\bf x}$ for every ${\bf x}\in{\mathbb R}^3\setminus {\rm int}B$ and $f(B)=B$. 
\end{enumerate}

\end{Proposition}

\vskip 3mm

\noindent{\bf Proof.} Up to conjugation we may suppose without loss of generality that $B$ is equal to the unit $3$-ball ${\mathbb B}^3$. Let $N_1$ be a regular neighbourhood of $V_1$ contained in $B={\mathbb B}^3$. Then $N_1$ is a solid torus and $N_1\setminus{\rm int}V_1$ is homeomorphic to ${\mathbb S}^1\times{\mathbb S}^1\times[0,1]$ where ${\mathbb S}^1$ denotes the unit circle and $[0,1]$ denotes the unit interval. Let $\varphi:{\mathbb S}^1\times{\mathbb S}^1\times[0,1]\to N_1\setminus{\rm int}V_1$ be a homeomorphism. Let $g:V_1\to V_2$ be a homeomorphism sending $k_1$ to $k_2$. For ${\bf x}\in V_1$ we define $f({\bf x})=g({\bf x})$. For ${\bf x}\in{\mathbb R}^3\setminus{\rm int}N_1$ we define $f({\bf x})={\bf x}$. For ${\bf x}\in N_1\setminus{\rm int}V_1$ with ${\bf x}=\varphi(a,b,t)$, we define $f({\bf x})=(1-t)f(\varphi(a,b,0))+tf(\varphi(a,b,1))$. Then we have a continuous map $f:{\mathbb R}^3\to{\mathbb R}^3$. It is clear that $f$ satisfies the conditions (1) and (2) except the condition $f({\mathbb B}^3)={\mathbb B}^3$. Since ${\mathbb B}^3$ is convex we see that $f({\mathbb B}^3)$ is a subset of ${\mathbb B}^3$. Since $f$ maps $\partial{\mathbb B}^3$ identically onto $\partial{\mathbb B}^3$ and there exist no retraction from ${\mathbb B}^3$ to $\partial{\mathbb B}^3$ we see that the image $f({\mathbb B}^3)$ must be the whole ${\mathbb B}^3$. 
$\Box$

\vskip 3mm

\begin{Example}\label{solenoid}
{\rm
Let $V$ be a knotted solid torus in ${\mathbb R}^3$. Let $f:{\mathbb R}^3\to{\mathbb R}^3$ be a continuous map such that the restriction map $f|_V:V\to{\mathbb R}^3$ is injective, $f(V)\subset{\rm int}V$ and $f(V)$ is essential in $V$. Namely there exists no $3$-ball in $V$ containing $f(V)$. 
Suppose that $\partial V$ and $f(\partial V)$ are not parallel. 
Let $k$ be a core of $V$. Then we have an infinite sequence $k,f(k),f^2(k),\cdots$ of knots. By a well-known fact on satellite knot we see that no two of them belong to the same knot type. The existence of such $f$ is assured by Proposition \ref{neighbourhood}. 
}
\end{Example}

\vskip 3mm

\begin{figure}[htbp]
      \begin{center}
\scalebox{0.4}{\includegraphics*{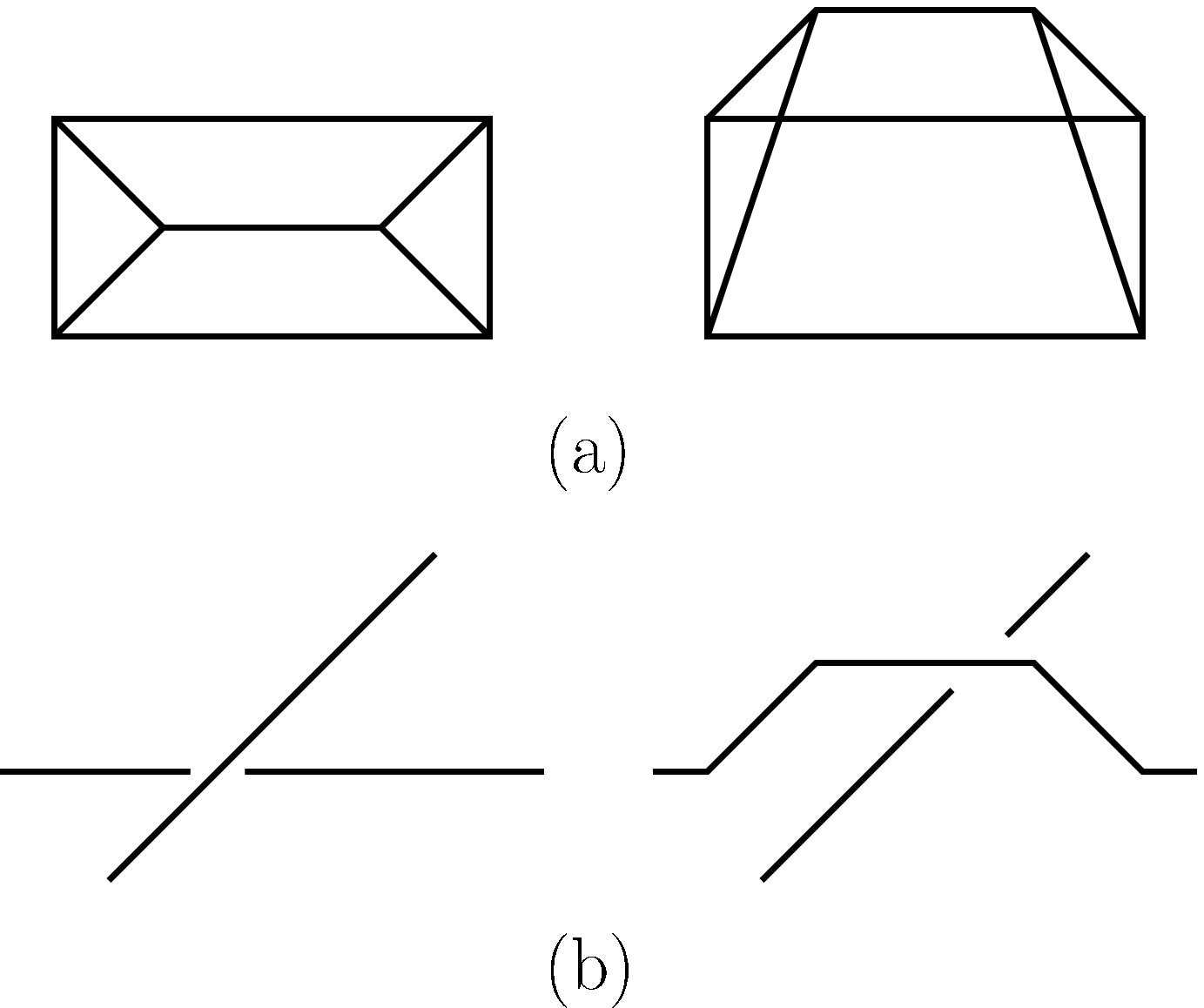}}
      \end{center}
   \caption{}
  \label{local-crossing-change}
\end{figure}

\begin{Example}\label{local crossing change}
{\rm
There exists a continuous map from ${\mathbb R}^3$ to ${\mathbb R}^3$ that is identical outside a small ball such that a crossing of a knot in the ball is changed by the map. To be more concrete, we give the following example. Let $f:{\mathbb R}^3\to{\mathbb R}^3$ be a continuous map defined as follows. 

(1) $f(x,y,z)=(x,y,z)$ if $|z|\geq1$ or $\sqrt{x^2+y^2}\geq2$,

(2) $f(x,y,z)=(x,y,3z+2)$ if $-1\leq z\leq0$ and $\sqrt{x^2+y^2}\leq1+|z|$,

(3) $f(x,y,z)=(x,y,-z+2)$ if $0\leq z\leq1$ and $\sqrt{x^2+y^2}\leq1+|z|$,

(4) $f(x,y,z)=(x,y,z-2\sqrt{x^2+y^2}+4)$ if $-1\leq z\leq1$ and $\sqrt{x^2+y^2}\geq1+|z|$. 

\noindent
See Figure \ref{local-crossing-change} (a) that illustrates $f$ on the $xz$-plane. 
Let $l_1$ be the $x$-axis $\{(x,0,0)|x\in{\mathbb R}\}$ and $l_2$ the line $\{(0,y,1)|y\in{\mathbb R}\}$. Then we see that $f(l_1)$ is above $f(l_2)=l_2$ at $(x,y)=(0,0)$ as illustrated in Figure \ref{local-crossing-change} (b). 
}
\end{Example}

Based on Proposition \ref{neighbourhood}, we can show the following theorem.  We denote the set of all positive integers by ${\mathbb N}$. 

\vskip 3mm

\begin{Theorem}\label{infinite sequence}
Let $\varphi:{\mathbb N}\to{\mathcal K}$ be a map. Then there is a continuous map $f:{\mathbb R}^3\to{\mathbb R}^3$ and an oriented polygonal knot $k$ in ${\mathbb R}^3$ that is iteratively injective with respect to $f$ such that $f^{n-1}(k)\in\varphi(n)$ for each $n\in{\mathbb N}$. 
\end{Theorem}

\vskip 3mm

\noindent{\bf Proof.} Let $B_n$ be a $3$-ball in ${\mathbb R}^3$ of radius $\frac{1}{3}$ centered at $(n,0,0)$ for each $n\in{\mathbb N}$. 
Let $k_n\in\varphi(n)$ be a polygonal knot contained in the interior of $B_n$ for each $n\in{\mathbb N}$. 
Let $g:{\mathbb R}^3\to{\mathbb R}^3$ be a map defined by $g(x,y,z)=(x+1,y,z)$. Then $k_n$ and $g^{-1}(k_{n+1})$ are oriented knots in a ball $B_n$. 
By Proposition \ref{neighbourhood} there exists a continuous map $g_n:{\mathbb R}^3\to{\mathbb R}^3$ such that $g_n(k_n)=g^{-1}(k_{n+1})$ and $g_n({\bf x})={\bf x}$ for every ${\bf x}\in{\mathbb R}^3\setminus{\rm int}B_n$. Let $f_n:B_n\to B_{n+1}$ be a continuous map defined by $f_n({\bf x})=g(g_n({\bf x}))$. Then we have $f_n(k_n)=k_{n+1}$ and $f_n({\bf x})=g({\bf x})$ for each ${\bf x}\in\partial(B_n)$. Now we define $f:{\mathbb R}^3\to{\mathbb R}^3$ by
\begin{enumerate}
\item[(1)] $f({\bf x})=g({\bf x})$ for ${\bf x}\in{\mathbb R}^3\setminus\bigcup_{n\in{\mathbb N}}B_n$, 
\item[(2)] $f({\bf x})=f_n({\bf x})$ for ${\bf x}\in B_n$ for each $n\in{\mathbb N}$. 
\end{enumerate}
Then we have $f(k_n)=k_{n+1}$ for each $n\in{\mathbb N}$. Let $k=k_1$. Then we have $f^{n-1}(k)=k_n\in\varphi(n)$ for each $n\in{\mathbb N}$ as desired. 
$\Box$

\vskip 3mm

Let $\mu$ be a positive real number. A {\it tent map} $t_\mu:{\mathbb R}\to{\mathbb R}$ is defined by
\[
t_\mu(z)=\left\{\begin{array}{cc}\mu z&(z\leq\frac{1}{2})\\-\mu z+\mu&(\frac{1}{2}\leq z).\end{array}\right.
\]
Tent maps are important examples in discrete dynamical systems theory. 
See for example \cite{Melo-Strien}. 
Therefore, as an example of a continuous map from ${\mathbb R}^3$ to ${\mathbb R}^3$ that generates a nontrivial discrete dynamical system, we consider a map $\widetilde{t_\mu}:{\mathbb R}^3\to {\mathbb R}^3$ defined by $\widetilde{t_\mu}(x,y,z)=(x,y,t_\mu(z))$. 
More generally we consider the following continuous maps. Let $f:{\mathbb R}\to{\mathbb R}$ be a continuous map. Then a continuous map $\tilde{f}:{\mathbb R}^3\to {\mathbb R}^3$ is defined by $\tilde{f}(x,y,z)=(x,y,f(z))$. 
We say that a map $f:{\mathbb R}\to{\mathbb R}$ is {\it piecewise linear} if ${\mathbb R}$ is a locally finite union of some closed intervals ${\mathbb R}=\bigcup_{\lambda\in\Lambda}I_\lambda$ such that the restriction map $f|_{I_\lambda}:I_\lambda\to {\mathbb R}$ is an Affine map for each $\lambda\in\Lambda$. By definition a piecewise linear map is continuous. 

\vskip 3mm

\begin{Proposition}\label{finiteness}
Let $f:{\mathbb R}\to{\mathbb R}$ be a piecewise linear map. Let $k$ be an oriented polygonal knot in ${\mathbb R}^3$ that is iteratively injective with respect to $\tilde{f}$. Then there exists a finite subset ${\mathcal J}$ of ${\mathcal K}$ such that the knot type of the oriented polygonal knot $(\tilde{f})^n(k)$ is an element of ${\mathcal J}$ for every non-negative integer $n$. 
\end{Proposition}

\vskip 3mm

For the proof of Proposition \ref{finiteness} we prepare the following Proposition \ref{same diagram}. It extends the well-known fact in knot theory that two knots with the same knot diagram are ambient isotopic. We call a knot diagram without over/under crossing information a {\it knot projection}. The condition (1) in Proposition \ref{same diagram} is a generalization of the condition that two knots have the same knot projection. Then the condition (2) in Proposition \ref{same diagram} is a generalization of the condition that two knots have the same over/under crossing information of their common knot projection. 
Let $p:{\mathbb R}^3\to{\mathbb R}^2$ be a natural projection defined by $p(x,y,z)=(x,y)$ and let $h:{\mathbb R}^3\to{\mathbb R}$ be a height function defined by $h(x,y,z)=z$. Let $k_1$ and $k_2$ be oriented polygonal knots in ${\mathbb R}^3$. We say that a homeomorphism $\varphi:k_1\to k_2$ is {\it polygonal} if there exists a subdivisions $k_1$ such that $\varphi$ maps each line segment of $k_1$ to a line segment of $k_2$. 

\vskip 3mm

\begin{Proposition}\label{same diagram}
Let $k_1$ and $k_2$ be oriented polygonal knots in ${\mathbb R}^3$. Suppose that there exists a polygonal homeomorphism $\varphi:k_1\to k_2$ that satisfies the following two conditions. 
\begin{enumerate}
\item[(1)] $p|_{k_1}=p|_{k_2}\circ\varphi$, 
\item[(2)] If ${\bf x},{\bf y}\in k_1$, ${\bf x}\neq{\bf y}$ and $p({\bf x})=p({\bf y})$, then $(h({\bf x})-h({\bf y}))(h(\varphi({\bf x}))-h(\varphi({\bf y})))>0$. 
\end{enumerate}
Then $k_1$ and $k_2$ are ambient isotopic. 

\end{Proposition}

\vskip 3mm

\noindent{\bf Proof.} For ${\bf x}\in k_1$ and $t\in[0,1]$ we define a point ${\bf x}_t$ by ${\bf x}_t=(1-t){\bf x}+t(\varphi({\bf x}))$. By the condition (1) we see that the points ${\bf x}$ and $\varphi({\bf x})$ are contained in a straight line parallel to the $z$-axis. 
Let $j_t=\{{\bf x}_t|{\bf x}\in k_1\}$. Then by the condition (2) we see that $j_t$ is a polygonal knot in ${\mathbb R}^3$. 
Since $j_0=k_1$ and $j_1=k_2$, $\{j_t|t\in[0,1]\}$ is an isotopy between $k_1$ to $k_2$. Clearly this isotopy has no non-locally flat points. Therefore it is extended to an ambient isotopy. See for example \cite[Theorem 3]{Hudson}. 
$\Box$

\vskip 3mm

It is well-known in knot theory that there are only finitely many knot types that share a common knot projection. 
We will see in the following proof of Proposition \ref{finiteness} that the knots $k, \tilde{f}(k), (\tilde{f})^2(k),\cdots$ have the same image under the natural projection $p:{\mathbb R}^3\to{\mathbb R}^2$. Note that the multiple points of them are not necessarily only finitely many transversal double points as in the usual knot projection. 
However it is still true that their knot types can vary only in a finite set as we will see in the following proof. 
By $[j]$ we denote the oriented knot type to which the oriented knot $j$ belongs. 
\vskip 3mm

\noindent{\bf Proof of Proposition \ref{finiteness}.} We set $k_n=(\tilde{f})^n(k)$. By the assumption we see that for $m\leq n$ the restriction map $(\tilde{f})^{n-m}|_{k_m}:k_m\to k_n$ is a polygonal homeomorphism satisfying the condition (1) in Proposition \ref{same diagram} and therefore $p(k_m)=p(k_n)$. If the condition (2) in Proposition \ref{same diagram} is also satisfied then $k_m$ and $k_n$ are ambient isotopic. 
We will show that there is a finite set $F$ of non-negative integers such that for every non-negative integer $n$ there exists $m\in F$ with $m\leq n$ so that $(\tilde{f})^{n-m}|_{k_m}:k_m\to k_n$ satisfies the condition (2) in Proposition \ref{same diagram}. Then the finite set ${\mathcal J}=\{[k_m]|m\in F\}$ is the desired set. 

Note that the set $S=\{(x,y)\in{\mathbb R}^2|{}^\#(p^{-1}(x,y)\cap k_n)\geq2\}$ is independent of the choice of a non-negative integer $n$. 
Let $k=l_1\cup\cdots\cup l_a$ be the line segments. In the following we pick up some essential cases and explain the idea of proof. 

First we consider the case that $p(l_i)$ is a point. Then $l_i$ is a line segment parallel to the $z$-axis. Then $(\tilde{f})^s(l_i)$ is also a line segment parallel to the $z$-axis for each non-negative integer $s$. 
By the assumption that the restriction map $(\tilde{f})^s|_k:k\to{\mathbb R}^3$ is injective, $(\tilde{f})^s$ maps $l_i$ homeomorphically onto $(\tilde{f})^s(l_i)$. Therefore $(\tilde{f})^{n-m}$ maps $(\tilde{f})^m(l_i)$ homeomorphically onto $(\tilde{f})^n(l_i)$. 
For $x, y\in l_i$ with $h(x)>h(y)$, 
$(h((\tilde{f})^{m}(x))-h((\tilde{f})^{m}(y)))(h((\tilde{f})^{n}(x))-h((\tilde{f})^{n}(y)))$ is positive if and only if both $h((\tilde{f})^{m}(x))-h((\tilde{f})^{m}(y))$ and $h((\tilde{f})^{n}(x))-h((\tilde{f})^{n}(y))$ are positive or both of them are negative. 
Note that for $z, w\in l_i$ with $h(z)>h(w)$, $h((\tilde{f})^{s}(z))-h((\tilde{f})^{s}(w))$ is positive if and only if $h((\tilde{f})^{s}(x))-h((\tilde{f})^{s}(y))$ is positive. 
After all, if both $h((\tilde{f})^{m}(x))-h((\tilde{f})^{m}(y))$ and $h((\tilde{f})^{n}(x))-h((\tilde{f})^{n}(y))$ are positive or both of them are negative then $(h(u)-h(v))(h((\tilde{f})^{n-m}(u))-h((\tilde{f})^{n-m}(v)))$ is positive for any $u,v\in (\tilde{f})^{m}(l_i)$ with $u\neq v$. 

Next we consider the case that there are line segments $l_i$ and $l_j$ with $i<j$ such that both $p(l_i)$ and $p(l_j)$ are line segments and $p(l_i)\cap p(l_i)$ is also a line segment. Note that $l_i$ and $l_j$ may or may not have a common end point. 
Let $x$ and $y$ be different points in $l_i$ and $l_j$ respectively with $p(x)=p(y)$. 
Let $z$ and $w$ be another pair of different points in $l_i$ and $l_j$ respectively with $p(z)=p(w)$. 
Since $(\tilde{f})^s$ maps $l_i\cup l_j$ injectively onto $(\tilde{f})^s(l_i\cup l_j)$ we see that $h((\tilde{f})^s(x))-h((\tilde{f})^s(y))$ is positive if and only if $h((\tilde{f})^s(z))-h((\tilde{f})^s(w))$ is positive. 
Therefore we have the following as in the previous case. 
If both $h((\tilde{f})^{m}(x))-h((\tilde{f})^{m}(y))$ and $h((\tilde{f})^{n}(x))-h((\tilde{f})^{n}(y))$ are positive or both of them are negative then $(h(u)-h(v))(h((\tilde{f})^{n-m}(u))-h((\tilde{f})^{n-m}(v)))$ is positive for any $u\in (\tilde{f})^{m}(l_i)$ and $v\in (\tilde{f})^{m}(l_j)$ with $u\neq v$ and $p(u)=p(v)$. 

As in these two cases we have the following. There are finitely many pairs of points $x_1,y_1,\cdots,x_a,y_a$ of $k$ with $p(x_i)=p(y_i)$ for each $i\in\{1,\cdots,a\}$ with the following properties. 
If both $h((\tilde{f})^{m}(x_i))-h((\tilde{f})^{m}(y_i))$ and $h((\tilde{f})^{n}(x_i))-h((\tilde{f})^{n}(y_i))$ are positive or both of them are negative for each $i\in\{1,\cdots,a\}$ 
then $(\tilde{f})^{n-m}|_{k_m}:k_m\to k_n$ satisfies the condition (2) in Proposition \ref{same diagram}. 
Thus we have the desired finite subset ${\mathcal J}$ of ${\mathcal K}$ with at most $2^a$ elements. $\Box$

\vskip 3mm

We say that a continuous map $f:{\mathbb R}\to{\mathbb R}$ is {\it switching} if it satisfies the following condition ($\ast$). 

\vskip 3mm

\noindent
($\ast$) For any map $\psi:{\mathbb N}\to\{-1,1\}$ there exists $x,y\in{\mathbb R}$ such that $f^{n-1}(x)\neq f^{n-1}(y)$ for each $n\in{\mathbb N}$ and $\displaystyle{\frac{f^{n-1}(x)-f^{n-1}(y)}{|f^{n-1}(x)-f^{n-1}(y)|}=\psi(n)}$ for each $n\in{\mathbb N}$. 

\vskip 3mm

For a switching continuous map $f:{\mathbb R}\to{\mathbb R}$ and a map $\psi:{\mathbb N}\to\{-1,1\}$, a pair $(x,y)$ of real numbers satisfying the condition ($\ast$) above is said to be a {\it realizing pair of $\psi$ with respect to $f$}. 

\vskip 3mm

\begin{Proposition}\label{switching map}
Let $\mu$ be a positive real number with $\mu>2$. Then a tent map $t_\mu:{\mathbb R}\to{\mathbb R}$ is switching. 
\end{Proposition}

\vskip 3mm

\noindent{\bf Proof.} It is well-known that there is a Cantor set $C_\mu\subset[0,1]$ with $t_\mu(C_\mu)=C_\mu$ such that for $x\in{\mathbb R}$, $\displaystyle{\lim_{n\to\infty}(t_\mu)^n(x)=-\infty}$ if and only if $x$ is not in $C_\mu$. 
Moreover the following is known. 
See for example \cite{Melo-Strien}. 
Let ${\mathcal S}$ be the set of all maps from ${\mathbb N}$ to $\{0,1\}$. For $x\in C_\mu$ we define $s(x)\in{\mathcal S}$ by 
\[
s(x)(n)=\left\{ \begin{array}{ll}
0 & \mbox{if}\ (t_\mu)^{n-1}(x)<\frac{1}{2}\\ \\
1 & \mbox{if}\ (t_\mu)^{n-1}(x)>\frac{1}{2}.\\
\end{array} \right.
\]
Then $s:C_\mu\to{\mathcal S}$ is a bijection. 
Let $d:{\mathcal S}\to{\mathcal S}$ be a map defined by $d(u)(n)=u(n+1)$. Then we have the following commutative diagram. 
\[\xymatrix{
C_\mu \ar[d]_{s} \ar[r]^{t_\mu|_{C_\mu}} & 
C_\mu \ar[d]^{s} & \\
{\mathcal S} \ar[r]_{d} & 
{\mathcal S} & 
}
\]
Let $\psi:{\mathbb N}\to\{-1,1\}$ be a map. Let $u:{\mathbb N}\to\{0,1\}$ be a map defined by
\[
u(n)=\left\{ \begin{array}{ll}
0 & \mbox{if}\ \psi(n)=-1\\ \\
1 & \mbox{if}\ \psi(n)=1\\
\end{array} \right.
\]
and let $v:{\mathbb N}\to\{0,1\}$ be a map defined by
\[
v(n)=\left\{ \begin{array}{ll}
0 & \mbox{if}\ \psi(n)=1\\ \\
1 & \mbox{if}\ \psi(n)=-1.\\
\end{array} \right.
\]
Let $x=s^{-1}(u)$ and $y=s^{-1}(v)$. Then for each natural number $n\in{\mathbb N}$ we have $(t_\mu)^{n-1}(x)<\frac{1}{2}<(t_\mu)^{n-1}(y)$ if $\psi(n)=-1$ and $(t_\mu)^{n-1}(x)>\frac{1}{2}>(t_\mu)^{n-1}(y)$ if $\psi(n)=1$. Therefore $\displaystyle{\frac{(t_\mu)^{n-1}(x)-(t_\mu)^{n-1}(y)}{|(t_\mu)^{n-1}(x)-(t_\mu)^{n-1}(y)|}=\psi(n)}$ for each $n\in{\mathbb N}$. Thus $t_\mu$ is switching. 
$\Box$

\vskip 3mm

\begin{Theorem}\label{tent map} Let $f:{\mathbb R}\to{\mathbb R}$ be a piecewise linear switching map. Let $\tilde{f}:{\mathbb R}^3\to {\mathbb R}^3$ be a continuous map defined by $\tilde{f}(x,y,z)=(x,y,f(z))$. Let ${\mathcal J}$ be any non-empty finite set of oriented tame knot types. Let $\varphi:{\mathbb N}\to{\mathcal J}$ be any map. Then there exists an oriented polygonal knot $k\subset{\mathbb R}^3$ that is iteratively injective with respect to $\tilde{f}$ such that $(\tilde{f})^{n-1}(k)\in\varphi(n)$ for each $n\in{\mathbb N}$.

\end{Theorem}

\vskip 3mm

\noindent{\bf Proof.} Let ${\mathcal J}=\{K_1,\cdots,K_m\}$ be a finite set of oriented knot types and $\varphi:{\mathbb N}\to{\mathcal J}$ a map. We may suppose without loss of generality that $\varphi(1)=K_1$. 
In the following we do not distinguish a knot diagram and its underlying projection so long as no confusion occurs. Here an underlying projection of a knot diagram is simply a subset of ${\mathbb R}^2$ and a knot diagram is an underlying projection together with over/under information at each crossing point. 
Let $D_i$ be an oriented knot diagram representing $K_i$ and ${\mathcal C}_i$ a set of crossings of $D_i$ such that changing all crossings in ${\mathcal C}_i$ will turn $D_i$ into a trivial knot diagram $D_i'$ for each $i\in\{1,\cdots,m\}$. Let ${\mathcal C}_i'$ be a set of crossings of $D_i'$ corresponding to ${\mathcal C}_i$ for each $i\in\{1,\cdots,m\}$. 
Let $D$ be a knot diagram obtained by a diagram-connected sum of $m$ knot diagrams $D_1,D_2',\cdots,D_m'$. 
Then $D$ is a diagram representing $K_1$. Let ${\mathcal C}$ be the set of all crossings of $D$. We may suppose that ${\mathcal C}_1,{\mathcal C}_2',\cdots,{\mathcal C}_m'$ are subsets of ${\mathcal C}$. Let ${\mathcal C}_0={\mathcal C}\setminus({\mathcal C}_1\cup{\mathcal C}_2'\cup\cdots\cup{\mathcal C}_m')$. 

Let $k$ be an oriented polygonal knot in ${\mathbb R}^3$ whose diagram is $D$. Namely $p(k)=D$ where $p(x,y,z)=(x,y)$ together with over/under crossing information. 
For a crossing $c\in{\mathcal C}$ let $o_c$ and $u_c$ be points in $k$ such that $p(o_c)=p(u_c)=c$ and $h(o_c)>h(u_c)$ where $h(x,y,z)=z$. 
Let $\psi_0:{\mathbb N}\to\{-1,1\}$ be a constant map defined by $\psi_0(n)=1$ for each $n\in{\mathbb N}$. 
For each $i\in\{1,\cdots,m\}$ let $\psi_i:{\mathbb N}\to\{-1,1\}$ be a map defined by $\psi_i(n)=1$ if $\varphi(n)=K_i$ and $\psi_i(n)=-1$ if $\varphi(n)\neq K_i$. 
By deforming $k$ by an ambient isotopy of ${\mathbb R}^3$ that preserves the $(x,y)$-coordinates if necessary, we may suppose that $k$ satisfies the following conditions. 
For each $i\in\{0,1\}$ and each crossing $c\in{\mathcal C}_i$, the pair $(h(o_c),h(u_c))$ of real numbers is a realizing pair of $\psi_i$ with respect to $f$, and 
for each $i\in\{2,\cdots,m\}$ and each crossing $c\in{\mathcal C}_i$, the pair $(h(u_c),h(o_c))$ of real numbers is a realizing pair of $\psi_i$ with respect to $f$. 
We now check that $(\tilde{f})^{n-1}(k)\in\varphi(n)$ for each $n\in{\mathbb N}$. 
For each crossing $c\in{\mathcal C}_0$ the condition that $(h(o_c),h(u_c))$ is a realizing pair of $\psi_0$ with respect to $f$ implies that the over/under crossing information of $(\tilde{f})^{n-1}(k)$ at $c$ is identical to that of $k$. 
For each crossing $c\in{\mathcal C}_1$ the condition that $(h(o_c),h(u_c))$ is a realizing pair of $\psi_1$ with respect to $f$ implies that the over/under crossing information of $(\tilde{f})^{n-1}(k)$ at $c$ is equal to that of $k$ if and only if $\varphi(n)=K_1$. 
For each $i\in\{2,\cdots,m\}$ and each crossing $c\in{\mathcal C}_i$ the condition that $(h(u_c),h(o_c))$ is a realizing pair of $\psi_i$ with respect to $f$ implies that the over/under crossing information of $(\tilde{f})^{n-1}(k)$ at $c$ is equal to that of $k$ if and only if $\varphi(n)\neq K_i$. 
After all we have that the diagram of the knot $(\tilde{f})^{n-1}(k)$ is a diagram-connected sum of $D_j$ with $\varphi(n)=K_j$ and $m-1$ trivial knot diagrams $D_i'$ with $i\in\{1,\cdots,m\}\setminus\{k\}$. Therefore knot $(\tilde{f})^{n-1}(k)$ is a representative of $\varphi(n)$ as desired. 
$\Box$

\vskip 3mm

\begin{Theorem}\label{universal map} There exists a continuous map $f:{\mathbb R}^3\to{\mathbb R}^3$ with the following properties. 
For any map $\varphi:{\mathbb N}\to{\mathcal K}$ there exists a tame knot $k\subset{\mathbb R}^3$ that is iteratively injective with respect to $f$ such that $f^{n-1}(k)$ is a tame knot in ${\mathbb R}^3$ and $f^{n-1}(k)\in\varphi(n)$ for each $n\in{\mathbb N}$. 

\end{Theorem}

\vskip 3mm

\begin{Remark}\label{remark}
{\rm The knot $k$ in the statement of Theorem \ref{universal map} may not be a polygonal knot nor a smooth knot. We do not know whether or not we can strengthen the condition on $k$ to be a polygonal knot or a smooth knot.}
\end{Remark}

\vskip 3mm

\noindent{\bf Proof of Theorem \ref{universal map}.} A closed subset $l$ of ${\mathbb R}^3$ is a {\it long knot} if it is abstractly homeomorphic to ${\mathbb R}$. By the one-point compactification of ${\mathbb R}^3$ we have a circle $l\cup\{\infty\}$ embedded in the $3$-sphere ${\mathbb R}^3\cup\{\infty\}$. Since the oriented knot types in the $3$-sphere are in one-to-one correspondence with those in the $3$-space we may think about the knot type of $l\cup\{\infty\}$ as an element of ${\mathcal K}$ unless it is not a tame knot. 
Let $t_3:{\mathbb R}\to{\mathbb R}$ be a tent map defined above and $\widetilde{t_3}:{\mathbb R}^3\to {\mathbb R}^3$ a continuous map defined by $\widetilde{t_3}(x,y,z)=(x,y,t_3(z))$. 
Let $X$ be the $x$-axis of ${\mathbb R}^3$. 
First we will construct a long knot $l\subset{\mathbb R}^3$ that is contained in a neighbourhood $N=\{(x,y,z)\in{\mathbb R}^3|y^2+z^2\leq9\}$ of $X$ so that $l$ satisfies the following condition. 
For each $n\in{\mathbb N}$ the knot type of $(\widetilde{t_3})^{n-1}(l)\cup\{\infty\}$ is equal to $\varphi(n)$. Actually $l$ is obtained from $X$ by countably many connected sum of possibly knotted arcs. Namely, for each $m\in{\mathbb N}$ we take a $3$-ball $B_m$ of radius $2$ centered at $(5m,0,0)$ and replace $X\cap B_m$ by a properly embedded arc $\alpha_m\subset B_m$ with $\partial\alpha_m\subset X$. Here $\alpha_1$ is knotted as the knot type $\varphi(1)$ and all other $\alpha_m$ with $m\geq2$ are unknotted. As in the proof of Theorem \ref{tent map} we can arrange each $\alpha_m$ so that $(\widetilde{t_3})^{n-1}(\alpha_m)$ is knotted as the knot type $\varphi(n)$ when $m=n$ and unknotted when $m\neq n$. 
Note that $(\widetilde{t_3})^{n-1}(\alpha_m)$ is not necessarily contained in $B_m$ but the knot type of it is still well-defined. 
Namely we may suppose that $(\widetilde{t_3})^{n-1}(\alpha_m)$ is a properly embedded arc in $[5m-2,5m+2]\times{\mathbb R}^2$ and then its knot type is naturally defined. 
It is easily checked by definition that the map $(\widetilde{t_3})^\ast:{\mathbb R}^3\cup\{\infty\}\to{\mathbb R}^3\cup\{\infty\}$ defined by $(\widetilde{t_3})^\ast({\bf x})=\widetilde{t_3}({\bf x})$ for ${\bf x}\in{\mathbb R}^3$ and $(\widetilde{t_3})^\ast(\infty)=\infty$ is continuous at $\infty$. 
Since ${\mathbb R}^3\cup\{\infty\}$ is homeomorphic to a $3$-sphere $({\mathbb R}^3\cup\{\infty\})\setminus\{(0,0,4)\}$ is again homeomorphic to the $3$-space ${\mathbb R}^3$. Let $h:({\mathbb R}^3\cup\{\infty\})\setminus\{(0,0,4)\}\to{\mathbb R}^3$ be a homeomorphism. By definition we see that the point $(0,0,4)$ is not in $\widetilde{t_3}({\mathbb R}^3)$. Then we have a continuous map $f:{\mathbb R}^3\to{\mathbb R}^3$ defined by $f({\bf x})=h((\widetilde{t_3})^\ast(h^{-1}({\bf x})))$. Let $k=h(l\cup\{\infty\})$. Then $k$ is an oriented knot in ${\mathbb R}^3$. The tameness of each knot $f^{n-1}(k)$ follows from the fact that the knot $(\widetilde{t_3})^{n-1}(l)\cup\{\infty\}$ is locally knotted only in a neighbourhood of $(\widetilde{t_3})^{n-1}(\alpha_n)$. This completes the proof. 
$\Box$

\vskip 3mm

\vskip 3mm

\section*{Acknowledgments} This work has been started on the opportunity of Professor Shin'ichi Suzuki's 70th birthday. 
The author would like to express his sincere gratitude to Professor Shin'ichi Suzuki for his constant guidance and encouragements.

{\normalsize
\end{document}